\documentclass[12pt,journal]{article}

\usepackage{cite}
\usepackage{amsmath,amssymb,amsthm}
\usepackage{bm}
\usepackage{dsfont}
\usepackage{graphicx}
\usepackage[caption=false,font=footnotesize]{subfig} 
\usepackage{epstopdf}   
\graphicspath{{figures/}}
\usepackage{amsmath}
\usepackage{amssymb}
\usepackage{mathrsfs}
\usepackage{turnstile}
\usepackage{xcolor}
\usepackage{tensor}
\usepackage{epstopdf} 
\usepackage{balance} 
\usepackage{todonotes}
\usepackage{ulem}

\usepackage{geometry}

\newcommand{\miu}{\leq}

\newcommand{\R}{\mathbb{R}}
\newcommand{\Z}{\mathbb{Z}}
\newcommand{\disp}{\displaystyle}
\newcommand{\be}{\begin{equation}}
\newcommand{\ee}{\end{equation}}
\newcommand{\phis}{\phi_{\sigma}}

\newcommand{\xx}{\underline{x}}

\newcommand{\kk}{\underline{k}}

\usepackage{hyperref} 
\hypersetup{
    colorlinks=true,
    linkcolor=blue,
    filecolor=magenta,      
    urlcolor=cyan,
    pdfpagemode=FullScreen,
    }

\usepackage{mathtools}

\title{Microwave Remote Sensing of Soil Moisture, Above Ground Biomass and Freeze-Thaw Dynamic:  Modeling and Empirical Approaches}
\author{L. Angeloni$^a$, D.D. Bloisi$^c$\thanks{Corresponding author: domenico.bloisi@unint.eu}, P. Burghignoli$^b$, D. Comite$^b$, D. Costarelli$^a$,\\ M. Piconi$^a$, A.R. Sambucini$^a$, A. Troiani$^a$, A. Veneri$^b$\\
\\
$^a$ Department of Mathematics and Computer Science,\\ University of Perugia,\\
$^b$ Department of Information Engineering, Electronics\\ and Telecommunications,
Sapienza University of Rome\\
$^c$ Faculty of Political Sciences and Psycho-Social Dynamics,\\
International University of Rome
}

\date{}

\begin{document}

\maketitle

\begin{abstract}
Human actions have accelerated changes in global temperature, precipitation patterns, and other critical Earth systems. Key markers of these changes can be linked to the dynamic of Essential Climate Variables (ECVs) and related quantities, such as Soil Moisture (SM), Above Ground Biomass (AGB), and Freeze-Thaw (FT) Dynamics. These variables are crucial for understanding global climate changes, hydrological and carbon cycles included. Monitoring these variables helps to validate climate models and inform policy decisions. Technologies like microwave remote sensing provide critical tools for monitoring the effects of human activities on these variables at a global scale. Other than proper tachenological developments, the study of ECVs requires suitable theoretical retrieval tools, which leads to the solutions of inverse problems. In this brief survey, we analyze and summarize the main retrieval techniques available in the literature for SM, AGB, and FT, performed on data collected with microwave remote sensing sensors. Such methods will be some of the fundamental algorithms that can find applications in the research activities of the interdisciplinary, curiosity-driven, project {\it REmote sensing daTa INversion with multivariate functional modeling for essential climAte variables characterization (RETINA)}, recently funded by the European Union under the Italian National Recovery and Resilience Plan of NextGenerationEU, under the Italian Ministry of University and Research. The main goal of RETINA, in which three research units from three different italian universities are involved, is to create innovative techniques for analyzing data generated by the interaction of electromagnetic waves with the Earth's surface, applying theoretical insights to address real-world challenges.
\end{abstract}

\section{Introduction}
Human activities such as burning fossil fuels, deforestation, agriculture, and industrial processes, are responsible for releasing significant amounts of carbon dioxide ($CO_{2}$), methane ($CH_{4}$), and other greenhouse gases, driving the rapid modifications the Earth's climate is experiencing.

These actions have caused long-term changes in temperature, precipitation patterns, and other Earth system dynamics.
Key markers of these changes include variations in Essential Climate Variables (ECVs) such as Soil Moisture (SM), Above Ground Biomass (AGB), and Freeze-Thaw (FT) Dynamics
(\url{https://gcos.wmo.int/en/essential-}\\
\url{climate-variables/about}). 
SM dynamics, a crucial part of the global hydrological cycle, are impacted by human-induced climate changes. These dynamics affect water availability, agricultural productivity, and even natural disaster patterns (e.g., droughts and floods).
Large-scale deforestation reduces the Earth's capacity to sequester carbon, directly impacting carbon cycles and exacerbating climate change.
AGB is a critical measure to monitor carbon stocks and understand the effects of deforestation.
Finally, the FT dynamics in polar and boreal regions are influenced by human activity emissions, accelerating permafrost melting and methane release, both of which have significant implications for climate change.

A central role to track these changes is played by Earth Observation technologies such as Microwave Remote (MW) Sensing (RS).
The monitoring of variables, like SM, AGB and FT, helps to validate climate models, understand the feedback mechanisms between human activity and environmental responses, and inform policies aimed at mitigating human impacts on the climate.
MW RS utilizes active sensors (e.g., radar) and passive sensors (e.g., radiometers) for continuous monitoring, irrespective of weather or lighting conditions. Depending on the platform (i.e., airborne or spacebone) data collection on regional and global scales is possible \cite{ulaby}. 

The retrieval techniques on MW RS data are based on theoretical \cite{theoretical}, semi-empirical \cite{semi-empirical}, and empirical \cite{empirical} models, increasingly enhanced by modern machine learning (ML) techniques utilizing neural networks (NNs). 
Theoretical models leverage classical electrodynamics theories, such as. e.g., scattering theory from rough surfaces \cite{ulaby}, Radiative Transfer Theory (RTT) \cite{RTT}, Foldy-Distorted Born Approximation (DBA) \cite{FDBA}.
Semi-empirical models, such as the Water Cloud Model (WCM) \cite{WCM}, links microwave signals to soil and vegetation parameters. Empirical models, instead, use direct relationships between observed signals and ECVs, they are therefore driven by the physical observable characterized by the microwave sensor.
ML is applied to handle complex datasets for more accurate predictions. However, the datasets required for the training are frequently very large; therefore, significant effort is needed for data annotation. In addition, data and ancillary data are often not continuously available due to acquisition methods, the type of sensor, the spatial and temporal resolition, as well as some practical conditions (e.g., the satellite orbit and the presence of disturbances like clouds).

\begin{figure}[h]
\includegraphics[width=8cm]{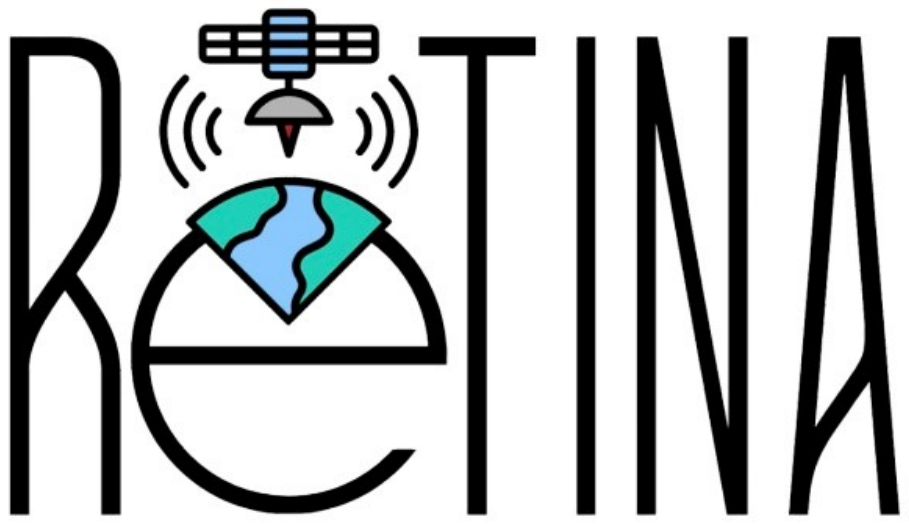}
\centering
\caption{The official logo of the \textit{RETINA} project.} \label{fig1}
\end{figure}

Very recently, a new research project thart aims to propose new methods for the retrieval of the ECVs has been funded, mixing both deterministic and non-deterministic procedures. This endeavor, called \textit{RETINA} (see Fig. \ref{fig1}, in which the official logo of the \textit{RETINA} project is provided), is
funded by the European Union in the frame of the Italian National Recovery and Resilience Plan (NRRP) of the NextGenerationEU program, under the Italian Ministry of University and Research. \textit{RETINA} proposes, for the first time, the application of direct and inverse analytical methods of Approximation Theory based on the theory of the so-called multivariate neural network (NN) operators (see, e.g., \cite{CS13,CCNP1,CP1}) for the modeling and estimation of SM, AGB, and FT, using data from space missions. The main idea is to combine techniques of analytic inversion (based on tools of functional analysis, such as series expansions \cite{howlett09}) as well as Bayesian approaches performed in conjuction with Monte Carlo methods. The main two (complementary) strategies of \textit{RETINA} can be summarized as follows: 
\begin{enumerate}
\item  Data Modeling with well-known Multivariate NN Operators. 
Through functional analysis techniques, theoretical inversion of these operators is achieved, resulting in an approximate analytical model for the target geophysical variables that is useful for their estimation. To address potential data disturbances, as well as take into account the uncertainty of the model, the NN operators will be extended to have the possibility of representing interval-valued fuzzy sets (IVFS), which allow for the representation of (uncertain variables) situations that are more coherent with real-world situations.

\item Bayesian Inversion with Monte Carlo Markov Chain (MCMC). 
Through MCMC techniques, Bayesian inversion complemented the NN operator approach. \textit{RETINA} targets the introduction of a specific type of Markov Chain, Probabilistic Cellular Automata (PCA), characterized by a parallel updating rule, which is expected to be particularly effective for retrieving multi-component physical quantities, such as matrix-formatted data.
\end{enumerate}

To set the state of the art of the retrieval techinques, a part of the RETINA project explores and summarizes the main algorithms available in the literature, focusing on the retrieval methods applied to the characterization of the bio-geophysical variables of interest. Providing an overview of these algorithms is the main motivation for this paper.

\section{Microwave remote sensing} \label{sec2}

MW RS is a powerful tool for monitoring bio-geophysical variables such as SM, AGB, and FT. This is given by capability to operate under all weather conditions, during day and night, and by the possibility of penetrating clouds, rain, and vegetation. Although often characterized by low spatial resolution (i.e., with respect to optical sensors), MW RS can generally offer relatively high temporal resolution with revisit time in the order of days, and up to hours in certain specific conditions (e.g., in the presence of satellite constellations).

Microwave signals are sensitive to the dielectric properties of soil, water, and vegetation, and are capable of penetrating surface layers (in the order of centimeters, depending on the wavelength and on the structure of the media), making them suitable for observing both superficial and shallow phenomena. 

Active sensors, including Synthetic Aperture Radar (SAR) systems, emit microwave radiation and collect the reflected or scattered signals. Passive sensors, such as radiometers, measure naturally emitted thermal radiation. MW radiation wavelengths range from 1 mm to 1 meter and are divided into different frequency bands (e.g., L-band, C-band, X-band and P-band), which correspond to decreasing wavelength. In MW RS, L-band (15–30 cm wavelength) and C-band (4–8 cm wavelength) are commonly used thanks to their balance between penetration depth, spatial resolution, and the technological maturity of the system.
Longer wavelengths (such as L-band or P-band, the latter planned to be used in the future ESA's Biomass mission \cite{ESABiomass}) are preferred for dense vegetation and high-biomass areas, even though they may suffer signal saturation. On the other hand, shorter wavelengths (e.g., C-band) are effective for less dense biomass or canopy surface and agricultural observations. Missions like SMAP \cite{SMAP} and SMOS \cite{SMOS} use both active and passive microwave sensors to monitor SM globally. 


\section{RETINA dataset of microwave remote sensing data}

One of the main goals of the RETINA project is to release a freely accessible dataset of remote sensing data that can be used for training, testing, and benchmarking retrieval procedures for the considered ECVs.

\begin{figure}[h]
\includegraphics[width=9cm]{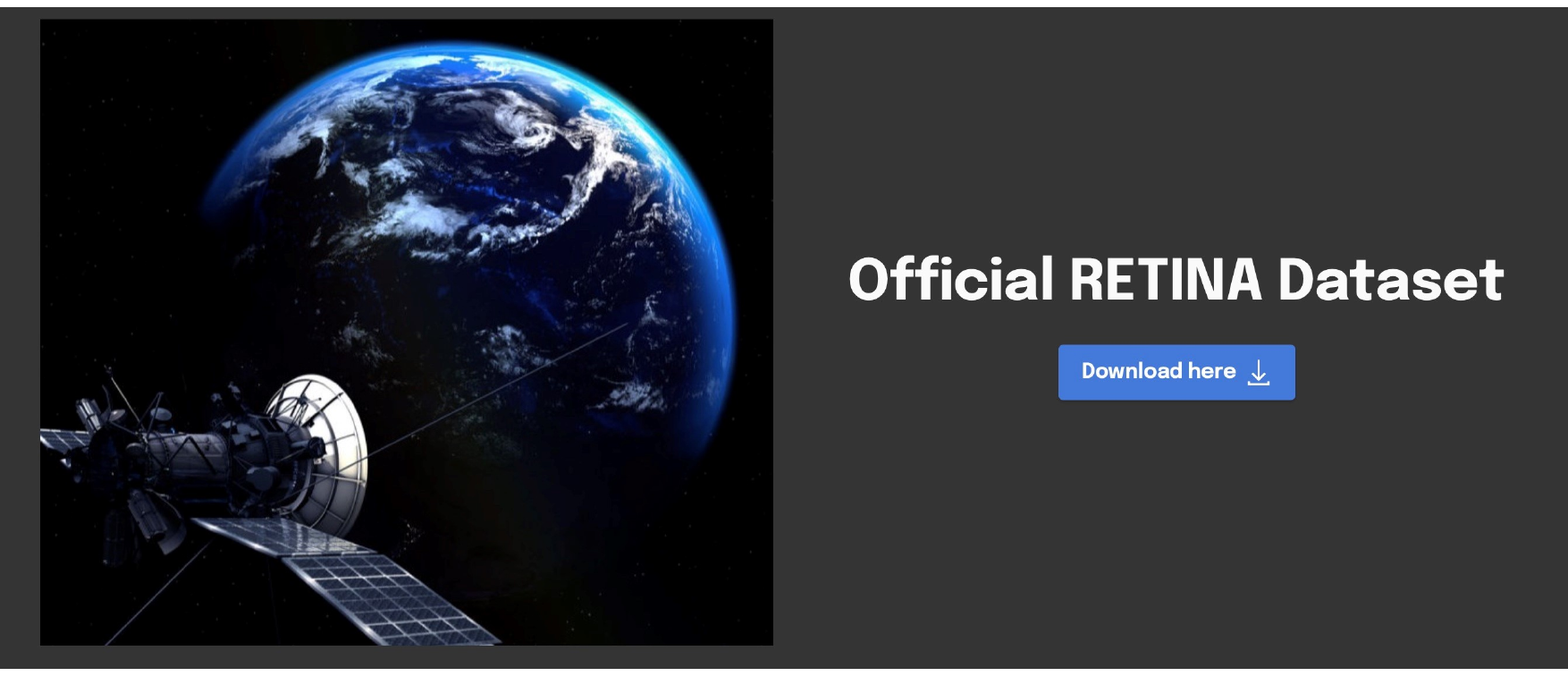}
\centering
\caption{The screen-shot of the RETINA dataset web page, available at the following link: \url{https://retina.sites.dmi.unipg.it/dataset.html}.} \label{Fig2}
\end{figure}

All the selected images have been collected in a dedicated open-access repository available from the RETINA website (see Fig. \ref{Fig2}) at the following link: \url{https://retina.sites.dmi.unipg.it/dataset.html}.

The source of all images (which can be downloaded for free, as stated in the copyright section of this paper) is the Sentinel-1 satellite constellation, the first of the ESA's Copernicus Program.
The source of all images (which can be downloaded for free, as stated in the copyright section of this paper) is the Sentinel-1 satellite constellation, the first of the ESA's Copernicus Program.
Images are available in GeoTIFF format with VV and VH bands.
To handle these images, it is possible to use any GIS software such as QGIS or SNAP.
Furthermore, a Python script for reading the images is 
available on the RETINA website (and provided below for convenience). To use the script, the packages "numpy", "rasterio" and "matplotlib" are required.
The script allows the user to clip the pixel intensities between user-chosen minimum and maximum threshold values to customize how images are displayed.

\begin{figure}[h]
\includegraphics[width=7cm]{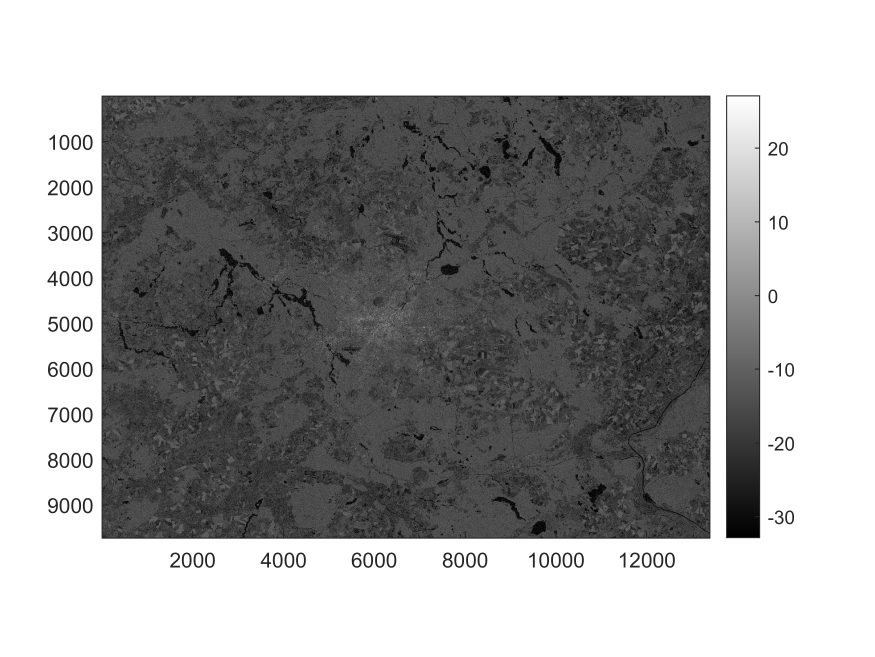}
\includegraphics[width=7cm]{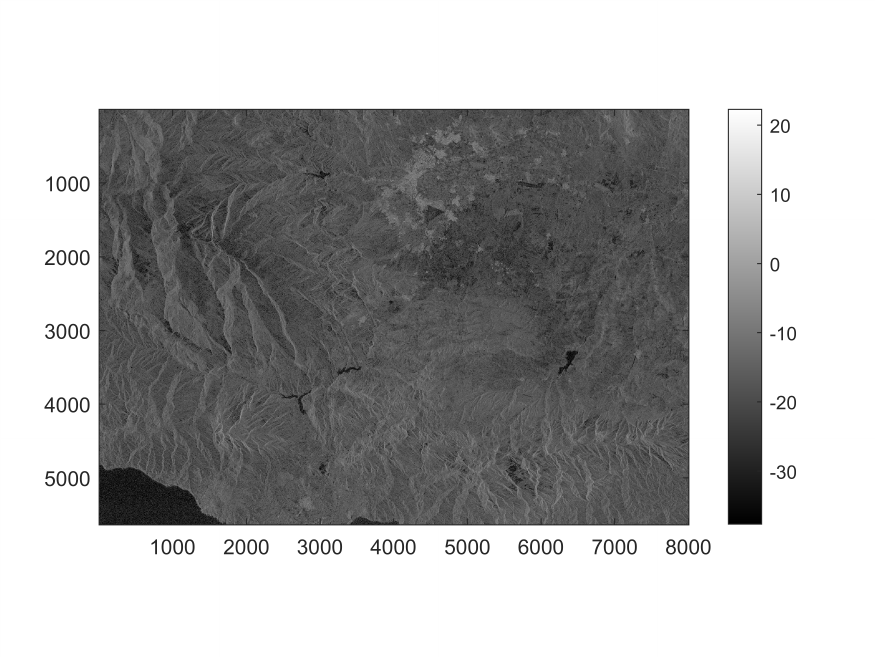}
\includegraphics[width=7cm]{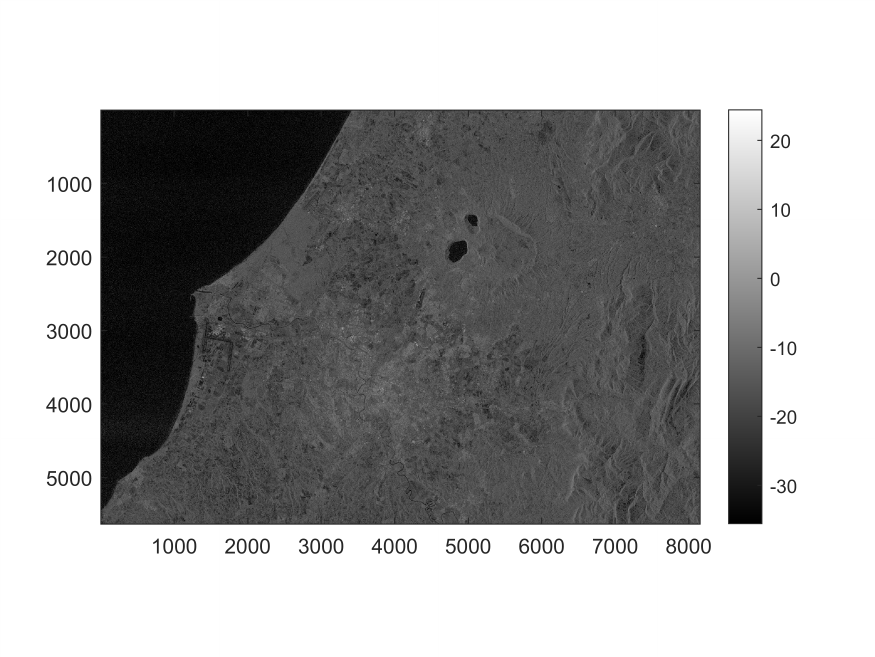}
\includegraphics[width=7cm]{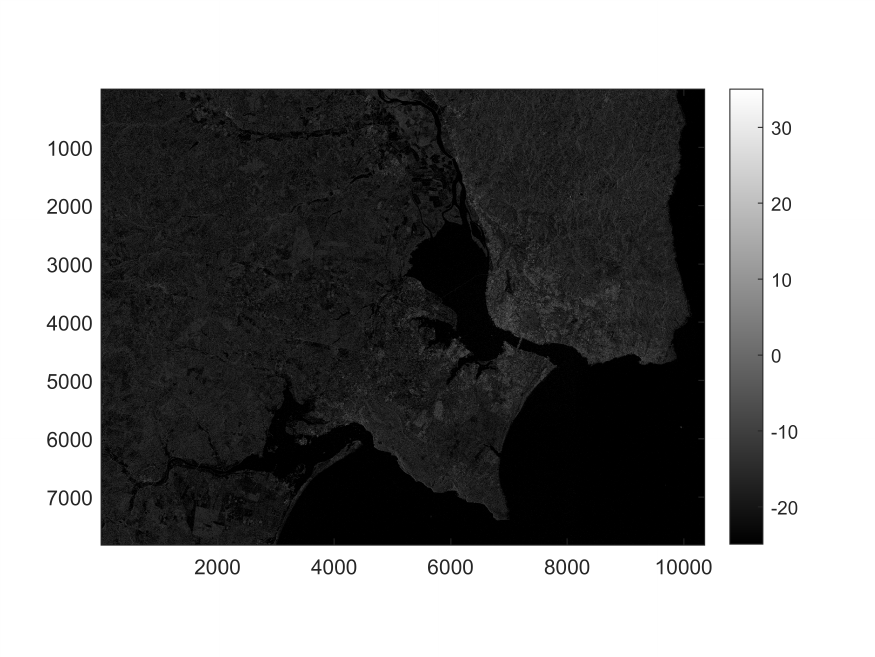}
\centering
\caption{RS images extracted from the RETINA dataset. On the top (from left to right): Berlin (Germany) and Granada (Spain) areas. On the bottom: (from left to right): Rome (Italy) and Lisbon (Portugal) areas.} \label{fig3}
\end{figure}

\begin{verbatim}
Python Code:
# Print the number of bands
print(f'Number of bands: {num_bands}')

# Optionally, print additional information
print(f'Width: {dataset.width}')
print(f'Height: {dataset.height}')
print(f'CRS: {dataset.crs}') # Coordinate Reference System (CRS)
print(f'Bounds: {dataset.bounds}')

# BAND STATISTICS
with rasterio.open(file_path) as dataset:
band1 = dataset.read(1)
band2 = dataset.read(2)

def print_stats(band, band_name):
print(f"{band_name} statistics:")
print(f" Min: {np.min(band)}")
print(f" Max: {np.max(band)}")
print(f" Mean: {np.mean(band)}")
print(f" Std: {np.std(band)}")

print_stats(band1, "Band 1")
print_stats(band2, "Band 2")

# Function to clip pixel values
def clip_array(array, min_percentile=2, max_percentile=98):
min_val = np.percentile(array, min_percentile)
max_val = np.percentile(array, max_percentile)
array_clipped = np.clip(array, min_val, max_val) 
# We keep the values of the array between min and max percentile
return array_clipped

# Clipping the bands
band1_clipped = clip_array(band1)
band2_clipped = clip_array(band2)

# Normalizing the clipped bands
band1_normalized = (band1_clipped - np.min(band1_clipped)) / 
(np.max(band1_clipped) - np.min(band1_clipped))
band2_normalized = (band2_clipped - np.min(band2_clipped)) / 
(np.max(band2_clipped) - np.min(band2_clipped))

plt.figure(figsize=(10, 10))
plt.title("Band 1 (Clipped and Normalized)")
plt.imshow(band1_normalized, cmap='gray')
plt.show()

plt.figure(figsize=(10, 10))
plt.title("Band 2 (Clipped and Normalized)")
plt.imshow(band2_normalized, cmap='gray')
plt.show()
\end{verbatim}
Also, a MATLAB script is available on the {\em RETINA} website. Examples of RS images that can be found in the RETINA dataset are provided in Fig. \ref{fig3}.

Below, we provide a short review of the main algorithms available in the literature for each one of the target ECVs considered in {\it RETINA}.

\section{Vegetation Biomass Retrival: an Overwiew} \label{sec4}

Retrieving vegetation biomass, specifically AGB, is crucial for understanding carbon cycles, forest dynamics, and the role of vegetation in climate change mitigation. MW RS, particularly radar systems, has emerged as a key method for estimating AGB due to its ability to penetrate vegetation layers and provide detailed structural information \cite{AGBMW}.

Longer Wavelengths (L-band, P-band) penetrate deeper into vegetation, interacting with trunks and larger branches, making them ideal for high-biomass regions. Shorter Wavelengths (e.g.,C-band) are sensitive to canopy features like leaves and smaller branches.

Techniques for Biomass Retrieval include:
\begin{itemize}
    \item  Theoretical approaches, such as the Radiative Transfer Theory (RTT) and the Wave Theory.
    RTT focuses on the principles of energy conservation and the transfer of energy through disordered media (i.e., layers of vegetation, see, e.g., \cite{kurum-2011,chen-2024}). A notable example is the Michigan Microwave Canopy Scattering model (MIMICS) \cite{Ulaby1990}, which divides the vegetation into crown, trunk, and ground components, and solves the vector radiative transfer equation (VRTE) iteratively. The Wave Theory, on the other hand, approximates solutions to Maxwell's equations to describe scattered electromagnetic fields. An example of this class is the Distorted Born Approximation (DBA), which uses electromagnetic wave theories of scattering to simulate how microwaves interact with vegetation structures. 
    DBA represents the vegetation layer as a random collection of individual scatterers (see, e.g., \cite{Lang2017}). These models are highly accurate but computationally intensive, also requiring a large collection of input data.
    \item Semi-empirical, such as the Water Cloud Model (WCM), models are simplified models that relate the measured radar backscattering coefficient to biomass using empirical relationships. WCM combines contributions from vegetation and ground backscatter (see, e.g., \cite{LI2018,NIJ2023}). These methods are adjusted for vegetation density and structure using parameters like canopy height and attenuation and are often used with regression techniques to estimate AGB from radar data. As the backscattering coefficient easily saturates at very high biomass densities, limiting the retrieval of AGB beyond a certain threshold, interferometric techniques are also often adopted as coherence appears to be effective to identify
    forested/nonforested areas and the height of the canopy \cite{Askne1997}. The main observable is the \textit{complex degree of coherence} (CDC), determining the potential of two electromagnetic signals to interfere \cite{Carminati2021}, which can be extracted from the data and compared with models (based on the WCM) accounting for stem volumes, tree height, and fill-factor, i.e., the fraction of ground covered by trees \cite{Santoro2002}.
    \item ML techniques utilize algorithms like NN and support vector machines to model complex relationships between radar signals and biomass (see \cite{ALI2015} for a complete review on this topic). They can integrate multi-source data (e.g., optical and microwave) to improve AGB retrieval accuracy.
\end{itemize}
SM, surface roughness, and vegetation water content can complicate signal interpretation. In general, retrieval models must account for different vegetation types, structures, and climates to ensure accuracy.

Current and upcoming radar missions for vegetation biomass studies are: 
ESA Biomass Mission (P-band SAR), designed specifically for global forest AGB mapping; NISAR (L- and S-band SAR), which targets forest structure and biomass dynamics; Sentinel-1 (C-band SAR) that provides data for biomass monitoring with limited penetration depth; GNSS Reflectometry, that is a more recent technique using reflected GPS signals (e.g., CYGNSS, TDS-1) to assess the AGB distribution at global scale (\cite{CYGNSS}).

\section{Freeze-Thaw Retrieval: An Overview}

FT retrieval focuses on monitoring seasonal transitions in the soil's thermal state, particularly between frozen and unfrozen conditions. These transitions are critical for understanding the water cycle, energy balance, and greenhouse gas dynamics, especially in high-latitude regions where permafrost melting can release significant amounts of methane \cite{Goulden_1998, Selvam_2016}.

MW RS has proven particularly well-suited for FT detection due to its high sensitivity to phase changes of water in soil, as the permittivity of water decreases dramatically between liquid and solid states \cite{ulaby}. 
Both active and passive MW remote sensing techniques are currently employed in the field, addressing the problem from different perspectives.
SAR and radiometer systems (e.g., Sentinel 1 and SMAP) provide high-resolution spatial data by capturing the dynamic temporal variations of the corresponding physical observable, which are associated with seasonal changes in FT states \cite{DERKSEN_2017}, \cite{KIMBALL_2004}. However, these systems are constrained by limited temporal resolution, typically on the order of days \cite{Luego_2019}.
%

FT retrieval methods include:
\begin{itemize}
    \item Theoretical models, mostly developed for passive radiometric measurements. This approach aims to estimate the brightness temperature from the physical properties of the medium by means of radiative transport theories. A notable example is the Helsinki University of Technology (HUT) model \cite{Lemmetyinen_2010}, which solves the scalar radiative transfer equation for multilayered systems.
    \item Empirical approaches, wherein thresholding is a commonly adopted method. For active measurements, the collected backscattering     coefficient is compared with its reference values for the thaw and freeze state, with the goal of determining the present water phase \cite{DERKSEN_2017}.  For GNSS systems the same method is adopted wherein the seasonal threshold algorithm now involves reference reflectivity values and measurements. Radiometric systems also employ a thresholding technique, but in this case, freeze/thaw discrimination is based on the difference between the vertical and horizontal polarizations of the brightness temperature.
    \item ML algorithms like Random Forest or NN process large datasets to detect patterns in FT transitions. These techniques are useful for integrating multi-sensor and auxiliary data (e.g., temperature and vegetation cover).
\end{itemize}
Problematics related to FT retrieval can include: vegetation and snow layers can obscure FT signals; rapid transitions may be missed without frequent observations; factors like soil roughness, composition, and moisture can complicate FT retrieval accuracy.
New systems like ESA Biomass will enhance FT monitoring capabilities.

\section{Soil Moisture Retrieval: An Overview} \label{sec6}
SM is a key bio-geophysical variable influencing global water cycles via evapotranspiration processes, exchange of heat between land and near-surface atmosphere, energy balance, and biochemical/carbon
cycles \cite{SENEVIRATNE2010125, SM}.
Recently, MW RS has emerged as a valuable tool for real-time SM monitoring due to its sensitivity to the soil-water ratio, which alters both the medium's emissivity and the backscattering properties of signals, opening the possibility for active and passive RS measurements. Key missions and instruments for SM study are: SMAP that combines SAR and radiometer data for global soil moisture monitoring; SMOS that uses passive L-band radiometer for large-scale SM and salinity observations; Sentinel-1 C-band SAR that provides high-resolution SM data, especially useful for agricultural applications.

Techniques for SM retrieval include:
\begin{itemize}
\item Theoretical models are particularly relevant in radar active measurements. Among them, the Integral Equation Method (IEM), based on electromagnetic scattering theory, has gained prominence over time and is now the most widely used \cite{fung_1992}. As active measurements are particularly affected by surface roughness conditions, the theory statistically accounts for the random variations of scattering surfaces, estimating the average values of scattered power as a function of surface roughness and the soil dielectric function.

\item Semi-empirical models, which are widely preferred due to their effectiveness and relatively simple implementation. For active systems, we distinguish two types of semi-empirical approaches. The first is based on the experimental calibration of the IEM, which corrects deviation of the theory from measurements correcting roughness effects \cite{baghdadi2002empirical}. The second exploits knowledge of scattering behavior in specific limiting cases, combined with experimental observations, to create ready-to-use formulas derived through data fitting. The most used frameworks are the Oh and the Dubois models \cite{Oh_1992, Dubois_1995}, both directly relating radar backscattering with volimetric soil moisture and roughness. 
For passive measurements, semi-empirical models are based on the scalar radiative transfer equation, solved at zeroth order and adjusted using experimental parameters to account for surface roughness, the mixing of different polarization components, and the influence of vegetation and atmospheric layers \cite{Wang_1981, KARTHIKEYAN2017106}. Semi-empirical models have also been used for GNSS-based SM retrieval, reveiling the potential of this novel technology to provide excellent results even at global scale \cite{Al-Khaldi_2019}.
\item ML and Data Assimilation: ML algorithms analyze complex, multi-sensor datasets to improve accuracy; data assimilation integrates satellite observations with hydrological models for comprehensive SM monitoring. For instance, Convolutional Neural Networks (CNNs) are useful for processing SAR data\cite{hegazi2021convolutional}, as they can extract spatial features from radar images, providing enhanced accuracy. More complex NN architecture can describe spatial features along with also temporal evolution.
\end{itemize}

Dense vegetation and uneven surfaces can obscure SM signals, so ground-truth data are often requested to ensure accuracy in various terrains and climates.

\section{Future developments: Neural Network Operators and Bayesian Inversion} \label{sec7}

NNs have become highly popular due to their utility across numerous fields, including Artificial Intelligence (AI), ML, and Approximation Theory (ATh). Within the \textit{RETINA} project, NNs will be applied in the context of ATh to develop approximate analytical models and their inversions for specific ECVs. In relation to the theory of NNs, in \cite{CS13}, the authors explored the functional properties of the NN operators. This research highlights the potential of these operators in modeling general two-dimensional structures, such as SAR satellite images or, more in general, RS data. 

For the sake of completeness, we recall the definition of such operators in both their classical and Kantorovich form. The multivariate discrete NN operators can be defined as follows:
\begin{equation}
     F^d_n(f,\xx)\ :=\ \frac{\disp \sum_{\kk \in {\cal J}_n} f\left({\kk \over n} \right)
        \Psi_{\sigma}(n \xx - \kk)}
      {\disp \sum_{\underline{j} \in {\cal J}_n}               \Psi_{\sigma}(n\xx - \kk)},  \quad \xx \in Q^d:=[a_1,b_1]\times \dots \times[a_d,b_d], 
\end{equation}
where the function:
\be
\Psi_{\sigma}(\xx) := \phi_{\sigma}(x_1) \cdot \phi_{\sigma}(x_2) 
           \cdots \phi_{\sigma}(x_d), 
                  \quad \quad    \xx := (x_1, ..., x_d) \in \R^s.
\ee 
is the multivariate (tensor-product) density function defined by means of suitable sigmoidal function $\sigma: \R \to \R$ (\cite{CY}), the set of indexes 
$$
{\cal J}_n := \left\{\kk \in \Z^d:\ \lfloor n a_i \rfloor \miu k_i \miu \lceil n b_i  \rceil  \right\}
$$
and
$$
\phis(x)\, :=\, \frac12\, \left[ \sigma(x+1)-\sigma(x-1) \right], \quad x \in \R.
$$.
We recall that, by Cybenko's definition in \cite{CY}, $\sigma: \R \to \R$ is called a {\em sigmoidal function} if $\lim_{x \to -\infty}\sigma(x)=0$ and $\lim_{x \to +\infty}\sigma(x)=1$.

While, the Kantorovich NN operators are:
\begin{equation}
     K^d_n(f,\xx)\ :=\ \frac{\disp \sum_{\kk \in {\cal J}_n} \left[ n^d \int_{R^n_{\kk}} f
           \left(\underline{u}\right) \, d\underline{u}\right]
        \Psi_{\sigma}(n \xx - \kk)}
      {\disp \sum_{\underline{j} \in {\cal J}_n}               \Psi_{\sigma}(n\xx - \kk)},  \quad \xx \in Q^d, 
\end{equation}
where:
\be
R^n_{\kk} := 
     \left[ \frac{k_1}{n}, \frac{k_1+1}{n} \right] \times  \cdots \times 
   \left[ \frac{k_d}{n}, \frac{k_d+1}{n} \right], 
\ee
are suitable multidimensional rectangles in which we will compute certain mean values of the considered function $f:Q^d\to \R$. 

We stress that, with respect to the classical (non-deterministic) theory of shallow and deep NNs, the NN operators are instead widely studied (see. e.g., \cite{Costarelli2022,CCNP1}) mathematical operators that are suitable to pursue a deterministic modeling approach, and also an enhancement/rescaling one. 

The task of inverting the above NN operators will be based on an analytical strategy, such as the possibility of exploiting Laurent's series, particularly when working with operators in Hilbert spaces, or methods of Approximation Theory. To address 2D data affected by measurement errors or other disturbances, interval-valued fuzzy sets (IVFS) have been proposed as a robust modeling tool (see \cite{alefeld2000}).

The NN operators approach will be enhanced through the integration of Bayesian inversion, which leverages advanced Monte Carlo Markov Chain (MCMC) techniques. These new techniques utilize a parallelized transition kernel to enable efficient sampling from the posterior distribution.
In the Bayesian framework, the variable to retrieve is treated as a random variable. With this approach, the retrieval procedure aims to determine the probability distribution of this variable given the observed data (posterior probability distribution). Then, the outcome of the retrieval procedure is the value that maximizes the posterior probability distribution. A standard way to estimate such probability distribution is to simulate the evolution of a Markov chain designed so that its stationary distribution is the posterior probability distribution of the variable to retrieve. Since the empirical distribution of the Markov chain converges to its stationary distribution, if the chain is run for a sufficiently long time, its empirical distribution is an estimate of the posterior probability distribution of the variables to retrieve. In this context, algorithms commonly used to simulate the Markov chain include the Metropolis algorithm, the Metropolis-Hastings algorithm, and the Gibbs sampler. 
If the Markov chain has multi-component states, such as in the case of 2D data, the previous algorithms sample the next state of the chain from a set of neighbors of the current state differing in only one component.
This strategy is referred to as single-flip dynamics \cite{hastings}.
An alternative to single-flip dynamics is the use of Probabilistic Cellular Automata (PCA) \cite{goldestein}. In PCA, all components of the state are updated simultaneously and independently at each step. This approach expands the set of neighbors to include the entire state space, resulting in higher motility and potentially faster convergence to equilibrium. 

More formally, a PCA, is a Markov chain $({X_n})_{n\in \mathbb{N}}$  defined on  $\mathcal{X} = \{1, \ldots, k\}^N$, where
$N$ is the number of components of the system, and whose transition matrix is such that
\begin{equation}
   \operatorname{\mathbb{P}}\{X_n = \tau | X_{n-1} = \sigma\} =
    \prod_{i=1}^N \operatorname{\mathbb{P}}\{(X_n)_{i} = \tau_{i} | X_{n-1} = \sigma\}.
\end{equation}
A transition matrix of this type is obtained by defining a function $H: \mathcal{X} \times \mathcal{X} \to \mathbb{R}$ of type
$H(\sigma, \tau) = -\sum_{1 \leq i \le N} h_{i}(\sigma)\tau_{i}$ and transition probabilities as
\begin{equation}
\operatorname{P(\sigma, \tau)} = \frac{e^{ -\beta H(\sigma, \tau) }}{Z_{\sigma}} = \prod_{1 \le i \leq N} \frac{e^{ \beta h_{i}(\sigma)\tau_{i} }}{(Z_{\sigma})_{i}}.
\end{equation}
where $\beta$ is a positive parameter called the inverse temperature and $Z_{\sigma}$ is a normalizing constant whose knowledge is not known to simulate the evolution of the chain.
Then, if $h(\cdot)$ satisfies certain suitable conditions, the stationary measure of the chain can be proven to be $\pi(\sigma) = \frac{\sum_{\tau}\operatorname{P}(\sigma, \tau)}{\sum_{\sigma, \tau}\operatorname{P}(\sigma, \tau)}$ (see articolo-da-scegliere).
The knowledge of the stationary measure of the chain allows to tune the algorithm so to favor the sampling of the more useful configurations (\cite{GibbsPCA, QUBO}).

  PCA's inherent parallelism offers significant computational advantages, particularly when leveraging massively parallel processing hardware such as GPUs or TPUs. These processors enable simultaneous updates of all components at each step, greatly enhancing efficiency.

\section{Conclusion} \label{sec8}

Climate change is one of the biggest challenges Mankind is called to face. The precise estimation and monitoring of the ECVs on a global scale is a fundamental step to describing and understanding the rapid changes the Earth's climate is experiencing and, consequently, to determine the more appropriate actions to mitigate the adverse effects of these changes.
  
Microwave remote sensing is a cornerstone of modern Earth observation, enabling critical insights into climate dynamics, environmental monitoring, and resource management at local, regional, and global scales.
This paper highlights advancements in microwave sensing technologies and their integration with ML to enhance the monitoring of Earth's bio-geophysical processes. 

Particular emphasis is given to those retrieval techniques that will be exploited to estimate the target ECVs in the context {\it RETINA} project.

\section*{Acknowledgments}

{\small The authors have been supported within the PRIN 2022 PNRR: "\textit{RETINA}: REmote sensing daTa INversion with multivariate functional modeling for essential climAte variables characterization", funded by the European Union under the Italian National Recovery and Resilience Plan (NRRP) of NextGenerationEU, under the Italian Ministry of University and Research (Project Code: P20229SH29, CUP: J53D23015950001).
}

\section*{Conflict of interest/Competing interests}

The authors declare no conflict of interest.

\section*{Availability of data and material and Code availability}

  The data made available in the framework of the {\it RETINA} project can be obtained
  from the "RETINA dataset" at \url{https://RETINA.sites.dmi.unipg.it/dataset.html}.
 {\it RETINA} Code and documentation are available through the website
  \url{https://RETINA.sites.dmi.unipg.itl}.
Permission is granted to use, copy, or modify the documentation for educational and research purposes without fee.

\section*{Copyright}
{\small The four images (Rome, Berlin, Lisbon, Granada) are available in the RETINA dataset via \url{https://retina.sites.dmi.unipg.it/dataset.html}. Permission to use, copy, or modify this dataset and its documentation for educational and research purposes only and without fee is granted, provided that this copyright notice and the original authors' names appear on all copies and supporting documentation. This dataset shall not be modified without first obtaining the permission of the authors. The authors make no representations about the suitability of this dataset for any purpose. It is provided "as is" without any express or implied warranty.
\\

\noindent In case of publishing results obtained using this dataset, please include explicit acknowledgments to the following website:
\\
\noindent \url{https://retina.sites.dmi.unipg.it/dataset.html}.


\end{document}